\documentclass[11pt]{article}
\usepackage{amsmath,amsthm,amsfonts}
\usepackage{color}
\numberwithin{equation}{section} \allowdisplaybreaks

\newtheorem{theorem}{\sc Theorem}[section]
\newtheorem{lemma}[theorem]{\sc Lemma}
\newtheorem{proposition}[theorem]{\sc Proposition}
\newtheorem{corollary}[theorem]{\sc Corollary}
\newtheorem{definition}[theorem]{\sc Definition}
\newtheorem{example}[theorem]{\sc Example}
\newtheorem{remark}[theorem]{\sc Remark}

\newcommand{\bet}{\begin{theorem}}
\newcommand{\eet}{\end{theorem}}
\newcommand{\blm}{\begin{lemma}}
\newcommand{\elm}{\end{lemma}}
\newcommand{\bprop}{\begin{proposition}}
\newcommand{\eprop}{\end{proposition}}
\newcommand{\bcor}{\begin{corollary}}
\newcommand{\ecor}{\end{corollary}}
\newcommand{\bdf}{\begin{definition}\rm}
\newcommand{\edf}{\end{definition}}
\newcommand{\bp}{\begin{proof}}
\newcommand{\ep}{\end{proof}}
\newcommand{\bex}{\begin{example}\rm}
\newcommand{\eex}{\end{example}}
\newcommand{\bremark}{\begin{remark}\rm}
\newcommand{\eremark}{\end{remark}}

\def\NN{{\mathbb N}}

\def\CC{{\mathbb C}}

\def\acc{{\rm acc}}

\newcommand{\codim}{\mathrm{codim \,}}

\begin{document}

\title {Pseudo-B-Fredholm operators, poles of the resolvent  and mean convergence in the Calkin Algebra   }

\author{ M. Berkani, S. \v{C}. \v{Z}ivkovi\'{c}-Zlatanovi\'{c}}

\date{}

\maketitle

\begin{abstract}

We define here a pseudo B-Fredholm operator as an operator such
that 0 is isolated in its essential spectrum, then  we prove that
an operator $T$ is  pseudo- B-Fredholm if and only if $T = R + F$
where $R$ is a Riesz operator and $F$ is a B-Fredholm operator
such that the commutator $[R,\, F]$ is compact. Moreover, we prove
that 0 is a pole of the resolvent of an operator $T$ in the Calkin
algebra if and only if $T= K+F$,  where $K$ is a power compact
operator and $F$ is a B-Fredholm operator, such that the
commutator $[K,\, F]$ is compact. As an application, we
characterize  the mean convergence in the Calkin algebra.

\end{abstract}

\renewcommand{\thefootnote}{}

\footnotetext{\hspace{-7pt}2010 {\em Mathematics Subject
Classification\/}:  primary 47A10, 47A53.
\baselineskip=18pt\newline\indent {\em Key words and phrases\/}:
 Pseudo-B-Fredholm, essential ascent, essential descent,  poles of the resolvent, Calkin algebra, mean convergence. }

\section{Introduction}
\noindent Let $X$ be an infinite dimensional Banach space  and let
$L(X)$ be the Banach algebra of bounded linear operators acting on
$X.$  The Calkin algebra over $X$ is the quotient algebra
$\mathcal{C}(X)=B(X)/K(X)$, where  $K(X)$ is the closed ideal of
compact operators on $X$. For $T\in L(X)$, let    $Ker(T)$ denote
the null-space and $R(T)$ the range of $T$. An operator $T \in
L(X)$ is Fredholm if $\text{dim} Ker(T)<\infty$)  and
$\text{codim} R(T)<\infty$.
 For $T\in L(X)$ {\em the Fredholm  spectrum},
 is defined by: $$ \sigma_F(T) = \{\lambda \in \CC: T-\lambda \text { is\ not\ a\ Fredholm\ operator}  \}$$

 Recall that the class of linear bounded B-Fredholm operators
were defined  in \cite{P7}. If $F_0(X)$ is the ideal of finite
rank operators in $ L(X)$ and
 $\pi: L(X)\longrightarrow  A $ is the canonical homomorphism,
where $ A= L(X)/F_0(X),$  it is well known by the Atkinson's
theorem \cite [Theorem 0.2.2, p.4]{BMW}, that $T \in L(X) $ is a
Fredholm operator if and only if its projection
 $\pi(T)$ in the   algebra $ A $ is invertible.  Similarly,   the
 following result  established  an Atkinson-type theorem for B-Fredholm
operators.

\begin {theorem}\cite [Theorem 3.4]{P10}\label{thm1} Let $T \in L(X).$ Then $T$ is a B-Fredholm operator
if and only if $\pi(T)$ is Drazin invertible in the algebra $
L(X)/F_0(X).$
\end{theorem}

\noindent  We conclude from the Atkinson's theorem and the
previous theorem, that invertibility in the algebra $ A=
L(X)/F_0(X)$ give rises to Fredholm operators,  while Drazin
invertibility in this algebra give rises to B-Fredholm operators.

Recall that an element  of   a unital  algebra $A$   is called
generalized Drazin invertible if $0$ is not an accumulation point
of its spectrum. Then it is natural to ask what are the properties
of those operators whose image under the canonical homorphism
$\pi: L(X)\longrightarrow A $ is generalized Drazin invertible?

Such operators will be called pseudo B-Fredholm operators, and
will be studied in the second section. The scheme
Fredholmness-Invertibility, B-Fredholmness- Drazin invertibility
is completed naturally by the couple pseudo-B-Fredholmness-
Generalized Drazin invertibility.

 In a recent works, among them \cite{BO}, \cite{TAJ} and
\cite{ZZ}, several authors studied pseudo-B-Fredholm operators as
the direct sum of a Fredholm operator and a quasi-nilpotent
operator. As will be seen, this definition is  a particular case
of our new definition, and by an example we prove the class of
pseudo-B-Fredholm operators we study here contains strictly the
class of operators studied by the authors cited above.

 In the main results of the second  section  we prove that the set $pBF(X)$ of
pseudo-B-Fredholm  operators in $L(X)$ is a regularity. Thus, by
usual properties of regularities, this implies that the
pseudo-B-Fredholm  spectrum $\sigma_{pBF}(T)$ satisfies the
spectral mapping theorem. Then we show that $T \in L(X)$ is a
pseudo-B-Fredholm operator  if and only if   $ T= R+F$ where $R$
is a Riesz  operator, $F$ is a  B-Fredholm operator such that the
commutator $[R, F]$ is compact or an inessential operator, and if
and only if $T$ is a compact perturbation of a direct sum of a
Fredholm operator and a Riesz operator.

 In the third section, we will answer successively the following three
 questions. The first one  is:  given a pseudo-B-Fredholm  operator $T,$
 when $0$ is a pole of its resolvent in the Calkin algebra $\mathcal{C}(X)$? We show
 that this holds if and only if $T= K+B,$  where $B$  is  a B-Fredholm
 operator and $K$ is a power compact one, such that the commutator $[K,\, B]$ is
compact.  The second question  is about
 the relation between the order of  $0$ as a pole of the resolvent of
 $\Pi(T)$ for a B-Fredholm operator $T$ and the essential ascent $a_e(T)$ and
the essential descent $ d_e(T)$ of $T,$ where $\Pi:
L(X)\longrightarrow L(X)/ K(X)$ is the canonical homomorphism. We
show that if $0$ is  a pole of order  $n,$ then $ n \leq a_e(T)=
d_e(T).$ Moreover we prove that $ n = a_e(T)= d_e(T)$ if and only
if $R(T^n)$ is  closed. The third question is: if $0$ is a pole of
the resolvent of $\Pi(T)$ of order $n$, when $T$ is then a
B-Fredholm operator with $ n = a_e(T)= d_e(T)$? The answer is that
this happens if and only if $ R(T^n)$ and $R(T^{n+1})$ are closed.
With the answer to those questions, we retrieve in particular some
similar results established in the case of Hilbert spaces in
\cite{BEL}.

As an application, we characterize mean ergodic convergence in the
Calkin algebra. Precisely, we show that the sequence $
(\Pi(M_n(T)))_n$ converges in the Calkin algebra if and only if
$\frac{ \mid\mid \Pi(T)^n\mid\mid}{n} \rightarrow 0$ as $
n\rightarrow \infty$ and there exists a power compact operator $K$
such that $ a_e(I-T+K)$ and $ d_e(I-T+K)$ are both finite and the
commutator $ [T, \, K] $ is compact, where $ M_n(T)= \frac{1+ T +
T^2+...T^n}{n}, n \in \NN$ and $ T \in L(X).$

We define now some tools that will be needed later. For $n\in\NN$
and $ T \in L(X),$ we set $c_n(T)=\dim R(T^n)/R(T^{n+1})$ and
$c_n^\prime(T)=\dim N(T^{n+1})/N(T^n)$. From \cite[Lemmas 3.1 and
3.2]{Kaashoek} it follows that $c_n(T)=\codim (R(T)+N(T^n))$ and
$c_n^\prime(T)=\dim (N(T)\cap R(T^n)).$   Obviously, the sequences
$(c_n(T))_n$ and $(c_n^\prime(T))_n$ are decreasing. The  {\it
descent} $\delta(T)$ and the {\it ascent} $ a(T) $ of $T$ are
defined by
 $ \delta(T)=\inf \{ n\in\NN:c_{n}(T)=0 \}=\inf  \{n\in\NN: R(T^n) = R(T^{n+1})\}$
and
 $ a(T)=\inf \{ n\in\NN:c^\prime_{n}(T)=0 \}= \inf \{n\in\NN : N(T^{n})=N(T^{n+1})\}$. We set formally $\inf\emptyset =\infty$.

The {\it essential  descent} $\delta_e(T)$ and the {\it essential
ascent} $ a_e(T) $ of $T$ are defined by
 $ \delta_e(T)=\inf \{ n\in\NN:c_{n}(T)<\infty \}$
and
 $ a_e(T)=\inf \{ n\in\NN:c^\prime_{n}(T)<\infty \}$.

Given a  Banach algebra $A$ and an element $a$ of $A,$  the left
multiplication operator $L_a: A \rightarrow A $ is defined by $
L_a(x)= ax,$ for all $ x \in A.$ It is well known that the
spectrum of $ a$ is equal to the spectrum of $L_a.$ We are
particularly interested in the case when $A$ is the Calkin algebra
and $a= \Pi(T)$ for $ T \in L(X).$

 The ascent and the descent of a Banach algebra element $ a \in
A$ are defined respectively as the ascent and the descent of the
operator $L_a.$

Now we give the definition of operators of topological uniform
descent, studied in \cite{GR}.

\bdf Let  $ T \in L(X) $ and let $  d \in {\bf N} $. Then $T$ has
a uniform descent for $ n \geq d$ if $R(T) + N(T^{n}) = R(T) +
N(T^{d}) $ for all $ n \geq d.$

 If in addition $ R(T) + N(T^{d}) $ is closed, then $T$ is said to
have a topological uniform descent for $ n \geq d$. \edf

 The {\it radical} of a unital Banach algebra $A$ is the set:
\begin{equation*}
\{d\in A:1-ad\ {\rm is\ invertible\ for\ all\ }a\in A\}=\{d\in
A:1-da \ {\rm is\ invertible\ for\ all\ }a\in A\} .\
\end{equation*}

 The set of all operators $A\in
L(X)$ satisfying $\Pi(A)\in{\rm Rad}(C(X))$, is the set of {\it
inessential} operators, denoted by $I(X)$.

For more details about those definitions, we refer the reader to
\cite{AI}.

\section{ Pseudo-B-Fredholm operators}

 \bdf \cite {KM} Let  $A$  be  an algebra over the field of complex numbers with a
unit $e$. A non-empty subset ${\bf{R}}$ of   $A$ is called a
regularity  if it satisfies the following conditions:
\begin{itemize}
 \item If $ a \in A$ and  $n\geq 1$ is an integer,  then
  $ a \in {\bf R} $ if and only if $a^n \in {\bf R}, $

 \item  If $ a,b,c,d  \in A $  are mutually commuting
  elements  satisfying  $ ac+bd = e, $ then $ ab\in {\bf R} $  if and
only if $a,b  \in {\bf R}. $

\end{itemize}

 \edf

Recall also that an element $ a \in A $ is said to be Drazin
invertible if there exists $b \in A $ such that $ bab=b, ab=ba$
and  $ aba-a$ is a nilpotent element in $A$.

\bdf An element  $a$  of a Banach algebra $A$ will be said to be
generalized Drazin invertible if there exists $b \in A $ such that
$ bab=b, ab=ba$ and  $ aba-a$ is a quasinilpotent element in $A$.
\edf

Koliha \cite{Koliha} proved that  $a\in A$ is generalized Drazin
invertible if and only if  there exists $ \epsilon >0,$ such that
for all $\lambda $ such that $ 0 < \mid \lambda \mid <  \epsilon,
$ the element  $a-\lambda e $ is invertible.

In the case of a general unital agebra, not necessarily a normed
algebra, we adopt this characterization as the definition of
generalized Drazin invertibility in such algebra. This is in
particular the case of the algebra $A= L(X)/ F_0(X).$

\bprop \label{same} Let $T \in L(X).$ Then $ \pi(T)$ is
generalized Drazin invertible  in the algebra $A=L(X)/F_0(X)$ if
and only if $\Pi(T)$ is generalized Drazin invertible in  the
Calkin algebra $\mathcal{C}(X).$ \eprop

\bp This  is a direct consequence of  the well known
characterization of Fredholm operators. \ep

\bdf \label{PBF} Let $T \in L(X).$ Then $T$ is said to be a
pseudo-B-Fredholm operator if $\pi(T)$ is generalized Drazin
invertible in the algebra $ L(X)/F_0(X).$ \edf

If $K \subset \mathbb{C}$, then $\acc \, K$ is the set of
accumulation points of $K$.

\bprop\label{bb1} Let $T \in L(X).$ Then $T$ is  a
pseudo-B-Fredholm operator if  and only if $0\notin\acc\,
\sigma_F(T)$. \eprop

\bp This  is a direct consequence of the Definition \ref{PBF}, the
characterisation of generalized Drazin invertible operators and
the characterization of Fredholm operators. \ep

 It is proved in
\cite{LUB} that the set of generalized Drazin invertible elements
in a unital  Banach algebra  is a regularity, from Proposition
\ref{same} we obtain immediately the following result.

\bet \label{} The set $pBF(X)$ of pseudo-B-Fredholm operators in
$L(X)$ is a regularity. \eet

Let  $\sigma_{pBF}(T)= \{ \lambda \in \mathbb{C} \mid T- \lambda I
\,\, \text {is not a pseudo-B-Fredholm operator }\} $  be the
spectrum generated by the regularity $pBF(X),$ for  $T \in L(X),$
then we have the following spectral mapping theorem.

 \bet  \label{Spectral-Theorem}  If   $f$  an analytic function in
a neighborhood  of the usual spectrum $\sigma(T)$  of \, an
operator $T$ in $L(X),$
 which is non-constant on any connected component of   $\sigma(T),$  then $
  f(\sigma_{{\bf pBF}}(T))=\sigma_{\bf pBF}(f(T)).$
 \eet

\bp This is a direct consequence of the properties of
regularities.\\ \ep

\bremark We  say that an operator $T \in L(X)$ is polynomially
Riesz if there exists a non-zero complex polynomial $P(z)$ such
that $P(T)$ is a Riesz operator. Every polynomially Riesz operator
in $L(X)$ is a pseudo-B-Fredholm operator.  Indeed if $ T$ is
polynomially Riesz, then  $P(T)$ is Riesz for a non-zero complex
polynomial $P(z) .$  As it is well known that the Fredholm
spectrum satisfies the spectral mapping theorem, then  we have $
P(\sigma_F(T))=\sigma_F(P(T))= \{0\} .$  Hence $\sigma_F(T)$ is
finite  because a polynomial has a finite set of roots. So it has
no  accumulation points  and from  Proposition \ref{bb1},  $T$ is
pseudo-B-Fredholm.

\eremark

For  $T \in L(X),$  we will say that  a subspace $M$
 of $X$  is {\em $T$-invariant} if $T(M) \subset M.$ We define
$T_{\mid M}:M \to M$ as $T_{\mid M}(x)=T(x), \, x \in M$.  If $M$
and $N$ are two closed $T$-invariant subspaces of $X$ such that
$X=M \oplus N$, we say that $T$ is {\em completely reduced} by the
pair $(M,N)$ and it is denoted by $(M,N) \in Red(T)$. In this case
we write $T=T_{\mid M} \oplus T_{\mid N}$ and say that $T$ is a
{\em direct sum} of $T_{\mid M}$ and $T_{\mid N}$.

 It is said that $T \in L(X)$ admits a  generalized
Kato decomposition, abbreviated as GKD,  if there exists $(M, N)
\in Red(T)$ such that $ T_{\mid M}$ is Kato and $T_{\mid N}$ is
quasinilpotent. Recall that an operator $T \in L(X)$ is {\em Kato}
if $R(T)$ is closed and $Ker(T) \subset R(T^n)$ for every $ n \in
\mathbb{N}$.

\bdf  $T \in L(X)$  is called a Riesz-Fredholm operator if
 there exists $(M, N) \in Red(T)$ such that $
T_{\mid M}$ is a Riesz operator and $T_{\mid N}$ is a Fredholm
operator.

\edf

 It is known that the sum of a Fredholm operator and a Riesz
operator   whose commutator is
 compact (or only  an inessential operator) is again a Fredholm operator. In the following theorem we show
 that an operator $T\in L(X)$ is  pseudo-B-Fredholm if and only  it is the  sum of a B-Fredholm operator and a Riesz operator
    whose commutator is a compact (an inessential) operator.

\bet \label{first} Let $T \in L(X).$ Then the following properties
are equivalent:

\begin{enumerate}

\item $T$  is a pseudo-B-Fredholm operator.

\item  $T$ is a compact perturbation of a Riesz-Fredholm operator.

\item $ T= R+B$ where $R$ is a Riesz  operator, $B$ is a
B-Fredholm operator such that the commutator $[R, B]$ is compact.

\item $ T= R+B$ where $R$ is a Riesz  operator, $B$ is a
B-Fredholm operator such that the commutator $[R, B]$ is an
inessential operator.
\end{enumerate}

\eet

\bp  $1)\Leftrightarrow 2)$ Suppose that $T$ is a
pseudo-B-Fredholm operator. If $T$ is Fredholm, then   the
statement  2 holds. Further, suppose that $T$ is not Fredholm,
then $0$ is
 an isolated point of $\sigma(\Pi(T)).$ Let $R \in C(X)$ be the spectral idempotent of $\Pi(T)$ corresponding to $\lambda=0,$
  then $R\ne 0$,  $\Pi(T)$ and $R$ commute, $\Pi(T)R$ is quasinilpotent and  $\Pi(T)+R$ is invertible according to  \cite[Theorem 3.1]{Koliha}.
   From \cite [Lemma 1]{BAR} we know that there exists
an idempotent $P \in L(X)$ such that $ \Pi(P)= R.$ Therefore,
$\Pi(TP)$ is quasinilpotent and $\Pi(T+P)$ is invertible.
 Since $\Pi(T)$ and $\Pi(P)$ commutes, we have that $\Pi(PTP)=\Pi(TP)$ and $\Pi((I-P)T(I-P))=\Pi(T(I-P))$. It follows that $PTP$
 is Riesz,
  $TP=PTP+K_1$, $T(I-P)=(I-P)T(I-P)+K_2$, where $K_1,K_2\in K(X)$, and
  so,
  \begin{eqnarray*}
    T &=& TP+ T(I-P)=PTP+(I-P)T(I-P)+K,
  \end{eqnarray*}
  \noindent where $K=K_1+K_2\in K(X)$.
    We have that $(R(P),R(I-P))\in Red(PTP)$,  $(R(P),R(I-P))\in Red((I-P)T(I-P))$,
  \begin{equation*}
    PTP = (PTP)_{\mid R(P)}\oplus (PTP)_{\mid R(I-P)}=(PTP)_{\mid R(P)}\oplus 0
  \end{equation*}
  and
  \begin{eqnarray*}
       (I-P)T(I-P) &=&((I-P)T(I-P)) _{\mid R(P)}\oplus ((I-P)T(I-P))_{\mid R(I-P)}\\&=&0\oplus ((I-P)T(I-P))_{\mid R(I-P)}.
  \end{eqnarray*}
  Therefore,
  \begin{equation}\label{first1}
    T=(PTP)_{\mid R(P)}\oplus ((I-P)T(I-P))_{\mid R(I-P)}+K.
  \end{equation}

  It's easily seen that $(PTP)_{\mid R(P)}$ is Riesz and further we prove that $((I-P)T(I-P))_{\mid R(I-P)}$ is Fredholm.
Since $\Pi(T+P)$ is invertible there exists $ S \in L(X)$ such
that $ \Pi(S)\Pi(T+P)=  \Pi(T+P)\Pi (S)= \Pi(I)$. As $\Pi(P)$ and
$\Pi(T+P)$ commute, then $\Pi(P)$ and $\Pi(S)$ commute and hence
\begin{eqnarray*}
  (I-P)(T+P)(I-P) (I-P) S(I-P)&=& I-P  +F_1,\\
  (I-P)S(I-P) (I-P) (T+P)(I-P)&=& I-P  +F_2,
\end{eqnarray*}
 where   $F_1$ and  $F_2$  are compact.
As $I-P$ is the identity on $R(I-P)$,  it follows   that
$((I-P)T(I-P))_{\mid R(I-P)}=((I-P)(T+P)(I-P))_{\mid R(I-P)}$ is a
Fredholm operator.
      According to \eqref{first1}, we see that  $T$ is a compact perturbation of
a Riesz-Fredholm operator.

 Conversely let $T=T_1\oplus T_2+K$ where $T_1$ is Riesz,
$T_2$ Fredholm and $K\in K(X)$.

It's clear  that $0$ is not an accumulation point of $
\sigma_F(T_1\oplus T_2)=\sigma_F(T)$ and according to Proposition
\ref{bb1},
 we get  that   $T$ is a pseudo-B-Fredholm operator.

$1)\Rightarrow 3)$ If $ T$ is a pseudo-B-Fedholm  operator, then
\begin{eqnarray*}
  T &=& (PTP)_{\mid R(P)}\oplus ((I-P)T(I-P))_{\mid R(I-P)}+K \\
  &=& [( (PTP)_{\mid R(P)}\oplus 0)+K]+[0\oplus ((I-P)T(I-P))_{\mid R(I-P)}],
\end{eqnarray*}
 where  $[( (PTP)_{\mid R(P)}\oplus 0)+K] $ is a Riesz operator and from \cite[Theorem 2.7]{P7}, $ [0\oplus ((I-P)T(I-P))_{\mid R(I-P)}] $ is a
B-Fredholm operator, here $P$ is the same idempotent as in  the
previous part of the proof. It is clear that the commutator of $ (
(PTP)_{\mid R(P)}\oplus 0)+K$ and $0\oplus ((I-P)T(I-P))_{\mid
R(I-P)} $ is compact.

$3)\Rightarrow 4)$ It follows from the inclusion $K(X)\subset
I(X)$.

$4)\Rightarrow 1)$ Let  $ T= R+B,$  where $R$ is a Riesz operator
and $B$ is a B-Fredholm  operator with $[R, B]$ is an inessential
operator. From \cite[Theorem 10.1]{ZDH} it follows that
\begin{equation}\label{b1}
 \sigma_{F}(T)=\sigma_{F}(B).
\end{equation}
  Since   $B$ is B-Fredholm, according to \cite[Remark A (iii)]{P12}  there
exists $ \epsilon > 0,$ such that if  $0<|\lambda| < \epsilon, $
we have that $ B- \lambda I$  is Fredholm which together with
\eqref{b1} gives that  $\lambda\notin \sigma_{F}(T)$. So
$0\notin\acc\,  \sigma_{F}(T)$ and thus $ T$ is a
pseudo-B-Fredholm operator by Proposition \ref{bb1}.
 \ep

\bigskip

 We mention that  Boasso considered in   \cite[Theorem 5.1]{BO}  isolated points of the spectrum of $ \Pi(T)$ for $ T \in L(X)$
 and he concluded  the equivalence
((1)$\Longleftrightarrow$(2))  by studying generalized Drazin
invertible elements in the range of a Banach algebra homomorphism
\cite[Theorem 3.2]{BO}, though our proof is more direct.

\begin{corollary}  Let $H$ be a Hilbert space and $ T \in L(H)$. Then $T$  is a pseudo-B-Fredholm operator if and only if
$ T= K+ Q + B,$ where $ K$ is compact, $Q$ quasi-nilpotent, $B$
B-Fredholm, with   $K$  and $[ Q, B]$ compact operators.
\end{corollary}
\bp In the case of  a Hilbert space, using the West decomposition
\cite{WEST} for a Riesz operator
 $R\in L(H)$ we have $ R= K + Q$ with $K$ compact and $Q$ quasi-nilpotent. Thus, according to Theorem \ref{first},
  $ T \in L(H)$  is a pseudo-B-Fredholm operator if and only if $ T= K+ Q + B,$ where $Q$ quasi-nilpotent, $B$  B-Fredholm, $K$
      and $[ Q, B]$ compact operators.
\ep

\bremark In the recent works \cite{BO}, \cite{TAJ} and \cite{ZZ},
the authors studied pseudo-B-Fredholm operators as the direct sum
of a Fredholm operator and a quasi-nilpotent one.  In
\cite[Theorem 3.4]{CZ} it is proved that
 \begin{eqnarray}
 & T \ {\rm is\  the\ direct\ sum\ of\ a\
Fredholm\ operator\ and\ a\ quasi-nilpotent\ one}\nonumber\\ &\Longleftrightarrow&\label{ds} \\& T\ {\rm admits\ a\ GKD\ and\ } 0 \not \in \acc \,
\sigma_{F}(T).\nonumber
 \end{eqnarray}
However there exists operators which are pseudo-B-Fredholm
operators in the sense of Definition \ref{PBF}, but do not have a
decomposition  as the direct sum of a Fredholm operator and a
quasi-nilpotent operator as seen by the following example.

Let $T$ be a compact   operator having infinite spectrum. Since $
\Pi(T)= 0,$ then $T$ is a pseudo-B-Fredholm operator   in the
sense of Definition \ref{PBF}.  We prove that $T$ cannot be
written as the direct sum of a Fredholm operator and a
quasi-nilpotent one. Assume the contrary. We observe first that $
T$ is not quasinilpotent, because it has non-zero spectrum. Also
$T$ is not Fredholm because is compact on the infinite dimensional
space $X$.

Assume that  there exists a pair $M, N$ of closed $T-invariant$
subspaces of $X$  such that $ T= T_1\oplus T_2 $ where
$T_1=T_{\mid M}$ is a quasi-nilpotent operator and $T_2=T_{\mid
N}$ is a Fredholm operator.  Since $\sigma(T)=\sigma(T_1)\cup
\sigma(T_2),$ it follows that $\sigma(T_2)$ is infinite. Therefore
$\ N$ is infinite-dimensional and  hence $\sigma_F(T_2)\ne
\emptyset$. As $T_2$ is Fredholm, it follows that
$0\notin\sigma_F(T_2)$.  Then  then there exists $ \lambda \neq 0
$ such that  $\lambda\in\sigma_F(T_2)$. But
$\sigma_F(T)=\sigma_F(T_1)\cup \sigma_F(T_2)$, $\sigma_F(T)=\{0\}$
and we get a contradiction. \eremark

 Our example shows that the condition that $T$ admits a GKD cannot be removed
  from the equivalence \eqref{ds}, in the other words the condition $ 0 \not \in \acc \,
\sigma_{F}(T)$ does not imply the condition that $T$ admits a GKD.

\section {Poles of the resolvent   in the Calkin algebra}

\bet \label{lifting} Let $T \in L(X)$. Then the following
properties are equivalent:

1- $0$ is a pole of the resolvent of $\Pi(T) $
 in the Calkin
algebra.

2- $T$ is the sum of a B-Fredholm operator $B$  and a power
compact operator $K$, such that the commutator $[K,\, B]$ is
compact.

3- There exists a  power compact operator $K$ such that $a_e(T+K)
$ and $d_e(T+K) $ are both finite, such that the commutator $[K,\,
T]$ is compact.

 \eet
\bp $1\Rightarrow 2)$ Suppose that $0$ is pole of the resolvent of
$\Pi(T) $ in the Calkin algebra. Let $R\in C(X)$ be the spectral
idempotent of $\Pi(T)$ corresponding to $\lambda=0$. Then $\Pi(T)$
and $R$ commute, $\Pi(T)R$ is nilpotent and $\Pi(T)+R$ is
invertible.
 From \cite [Lemma 1]{BAR} it follows  that there exists
an idempotent $P \in L(X)$ such that $ \Pi(P)= R.$  Since
$\Pi(TP)=\Pi(PTP)$, it follows that $PTP$ is a power compact
operator.
 As in the proof of Theorem \ref{first}, we get that there is $K'\in K(X)$ such that
  \begin{eqnarray*}
      T
     =(PTP)_{\mid R(P)}\oplus ((I-P)T(I-P))_{\mid R(I-P)}+K'.
  \end{eqnarray*}

  Set $ K= ((PTP)_{\mid R(P)}\oplus 0) + K'=PTP+K'$ and $ B= 0 \oplus ((I-P)T(I-P))_{\mid R(I-P)}.$
Then $K $ is clearly a power compact operator, $B$ is a B-Fredholm
operator by \cite[Theorem 2.7]{P7} and the commutator $[K,\, B]$
is compact.

$2)\Rightarrow 3) $ Assume that  $T$ is the sum of a B-Fredholm
operator $B$ and a power compact operator $K'$ such that the
commutator $[K',\, B]$ is compact.   Let $K= -K',$ then $ B= T+K,$
and  from \cite [Theorem 3.1]{P12}, $a_e(T+K) $ and $d_e(T+K) $
are both finite. Moreover, the commutator  $[K,\, T]$ is compact.

 $3)\Rightarrow  1)$ Assume that  there exists a  power compact operator $K$
such that $a_e(T+K) $ and $d_e(T+K) $ are both finite and the
commutator  $[K,\, T]$ is compact. Then from \cite [Theorem
3.1]{P12}, $ T+K$ is a B-Fredholm operator. Hence $ \Pi(T+K)$ is
Drazin invertible in the Calkin algebra. As the commutator  $[K,\,
T]$ is compact, then  $ \Pi(T) \Pi(K)= \Pi(K) \Pi(T).$ Since $ K$
is power compact, then $\Pi(K)$  is nilpotent.  From \cite[Theorem
3]{ZCCW} we know that Drazin invertibility is stable under
nilpotent commuting perturbations.  Thus  it follows that $\Pi(T)=
\Pi(T+K) -\Pi(K)$  is Drazin invertible in the Calkin algebra and
$0$ is pole of $\Pi(T).$

\ep

\bigskip

As a consequence, in the case of Hilbert spaces, we recover
\cite[Theorem 2.2]{BEL}.

\bcor\label{Hilbert-case}  Let $H$ be a Hilbert space and  $T \in
L(H)$. The following properties are equivalent:

1- $0$ is pole of the resolvent of $\Pi(T) $ in the Calkin
algebra.

2- There exist a  compact operator $K$ such that $ T+K$ is a
B-Fredholm operator.

 3-There exist a  compact operator $K$ such that $a_e(T+K)$ and $d_e(T+K) $ are both finite.

\ecor

\bp  $1)\Rightarrow 2)$ As in the proof of Theorem \ref{lifting},
there exist an idempotent $P \in L(H)$ and $K\in K(H)$ such that
$T=(PTP)_{\mid R(P)}\oplus ((I-P)T(I-P))_{\mid R(I-P)}+K$, where $
((I-P)T(I-P))_{\mid R(I-P)}$ is a Fredholm operator and  $PTP$ is
a power compact operator, which implies that $(PTP)_{\mid R(P)}$
is a power compact operator.   From \cite[Lemma 5]{HLY}, there
exists a nilpotent operator $N_1$ and a compact operator $K_1$
defined on $R(P)$  such that $ (PTP)_{\mid R(P)}= N_1 + K_1.$
Hence
$$ T= [N_1 \oplus ((I-P)T(I-P))_{\mid R(I-P)}] + [ (K_1 \oplus 0) +
K], $$  $(K_1 \oplus 0) + K$ is clearly a compact operator and by
\cite[Theorem 2.7]{P7}
 the operator $N_1 \oplus
((I-P)T(I-P))_{\mid R(I-P)}$ is a B-Fredholm operator.

 The proof of the other implications are similar to the
corresponding implications in Theorem \ref{lifting}. \ep

Recall that  from \cite[Theorem 5.3]{GRZ}, if $ T \in L(X)$ has
finite essential ascent $ a_e(T)$ and finite essential descent $
d_e(T),$ then they are equal. In the following result, for a
B-Fredholm operator $T$, we compare the order of $0$ as a pole of
the resolvent of $\Pi(T)$ and the common value of its essential
ascent and its essential descent.


\bet \label{order} Let $T \in L(X)$ be a B-Fredholm operator. Then
$\Pi(T)$ is Drazin invertible in the Calkin algebra and if $n$ is
the order of  $0$ as a pole of the resolvent of  $\Pi(T),$
 then $ n \leq a_e(T)= d_e(T).$ Moreover $ n = a_e(T)= d_e(T)$ if
and only if $R(T^n)$ is  closed

\eet

\bp Let $ d= a_e(T)= d_e(T)$  and assume that $d < n.$ Then from
\cite[Theorem 3.1]{P12} $R(T^d)$ is closed and the operator $T_d:
R(T^d) \rightarrow R(T^d)$ is a Fredholm operator. Thus there
exists a compact operator $K_{d}$  in $L(R(T^d)), $ an operator
$R_d$  in $L(R(T^d))$ such that $R_dT_d= I_d + K_{d}, $ where $
I_d$ is the identity of $L(R(T^d))$.  Thus  $ T^d-R_dT^{d+1}$ is a
compact operator.

Let  $ V \in L(X)$ such that $T^n V$ is compact. Then
$$ T^{n-1} V=(T^d-R_dT^{d+1}) T^{n-d-1}V +R_d T^{n}V$$ is a
compact operator. Hence $ a(\Pi(T)) \leq n-1 < n$ and this is a
contradiction. Thus $ n \leq  a_e(T)= d_e(T). $

If  $ n = a_e(T)= d_e(T),$ then  from \cite[Theorem 3.1]{P12}, $
R(T^n)$ is closed. Conversely if $ R(T^n)$ is closed, let $ d=
a_e(T)= d_e(T).$ Since $\Pi(T)$ is Drazin invertible in the Calkin
algebra, there exists $S \in L(X),$ such that the operators $
TS-ST, STS-S, T^nST-T^n$ are all compact operators. Let $ K=
ST^{n+1}- T^n,$ then $K$ is a compact operator and $ Ker(T) \cap
R(T^n) \subset R(K).$  As $Ker(T) \cap R(T^n)$ is closed, then
$Ker(T) \cap R(T^n)$ is finite dimensional. If $n < d,$ then $n
\leq d-1$ and $ Ker(T) \cap R(T^{d-1}) \subset Ker(T) \cap R(T^n)
.$  Hence $ Ker(T) \cap R(T^{d-1})$ is finite dimensional and
consequently $ a_e(T) \leq d-1,$ which is a contradiction. Hence $
n\geq d.$ As we know already that $ n \leq a_e(T)= d_e(T),$ then $
n = a_e(T)= d_e(T).$ \ep

\bremark Without the hypothesis of the closedness of the range
$R(T^n)$, Theorem \ref{order} may be false. For example let $K$ be
a nilpotent compact operator with infinite dimensional range, then
it is easily seen that the order of $0$ as a pole of $\Pi(K)$ is
equal to one, while $K$ has a finite essential ascent and descent
strictly greater than $1,$ because the range $R(K)$ of $K$ is not
closed.  \eremark

 \noindent In the
following result, we give a sufficient condition which implies the
equality of the order and the common value of the essential ascent
and descent. In the case of Hilbert spaces, this result was proved
in \cite[proposition 3.3]{BEL}. While in \cite[proposition
3.3]{BEL}, the proof is based on Sadovskii essential enlargment of
an operator \cite{SAD}, our proof is based directly on the
definition of the Drazin inverse.

\bet  Let $T \in L(X)$ be a B-Fredholm operator  with finite
essential ascent and finite essential descent equaling $d$. If
$d=0$ or $ R(T^{d-1})$ is closed, then $0$ is a pole of the
resolvent of $\Pi(T)$ of order $d.$

\eet

\bp  If $d= 0$ then $T$ is a Fredholm operator. So $\Pi(T)$ is
invertible in the Calkin  algebra and $0$   is a pole of the
resolvent of $\Pi(T)$ of order $d.$

 Assume now that $d> 0$ and $ R(T^{d-1})$ is closed. Since $T$
is a B-Fredholm operator, then
 $ \Pi(T)$ is Drazin invertible in the Calkin algebra. Let $ n=
a(\Pi(T))= d(\Pi(T)).$ Then there exists  $S \in L(X)$ such that
the operators $ TS-ST, STS-S, T^nST-T^n,$ are all compact
operators. Let $ K= ST^{n+1}- T^n,$  then $K$ is a compact
operator and $ Ker(T) \cap R(T^n) \subset R(K).$

We already know that  $  n\leq d.$ If $n < d,$ then $n \leq d-1$
and $ Ker(T) \cap R(T^{d-1}) \subset Ker(T) \cap R(T^n) \subset
R(K).$ As  $ Ker(T) \cap R(T^{d-1})$ is closed, then $ Ker(T) \cap
R(T^{d-1})$ is finite dimensional. Thus $ a_e(T) \leq d-1.$
Contradiction. Hence $ n= d.$

\ep

Now we give necessary and sufficient conditions to lift a Drazin
invertible element in the Calkin algebra as a B-Fredholm operator.
The sufficient condition in the case Hilbert spaces had been
proved in \cite[Theorem 3.2]{BEL}, where the proof is based also
on  Sadovskii essential enlargment of an operator \cite{SAD},
while we use here properties of B-Fredholm operators.

\bet \label{PolesCA} Let $T \in L(X)$  such that  $\Pi(T)$ is
Drazin invertible in the Calkin algebra and $0$ is a pole of the
resolvent of $\Pi(T)$ of order $n$.  Then $T$ is a B-Fredholm
operator with  $a_e(T)= d_e(T)= n$  if and only if
   $ R(T^n)$ and $R(T^{n+1})$ are closed.
\eet

\bp If  $T$ is a B-Fredholm operator with  $a_e(T)= d_e(T)= n, $
then from \cite[Theorem 3.1]{P12},     $ R(T^n)$ and $R(T^{n+1})$
are closed.

Conversely assume that  $ R(T^n)$ and $R(T^{n+1})$ are closed.
Since $0$ is a pole of the resolvent of  $\Pi(T)$  of order $n,$
then there exists an operator $S$  in $L(X),$ such that $ TS-ST,
STS-S, T^nST-T^n $ are all compact operators. Let $K=
ST^{n+1}-T^n,$ then $K$  is a compact operator. If $ y \in N(T)
\cap R(T^n),$ then $ y \in R(K).$ Since $ N(T) \cap R(T^n)$ is
 closed and $K$ is compact, then
$ N(T) \cap R(T^n)$  is of finite dimension. Let $T_n: R(T^n)
\rightarrow  R(T^n)$ be the operator induced by $ T$. Then $T_n$
is an upper semi-Fredholm operator and so, $T$ is a
semi-B-Fredholm operator. In particular and  from \cite{GR}  it is
an operator of
 topological uniform descent. Since  $0$ is a pole of the resolvent
  of  $\Pi(T)$  of order $n,$  if $|\lambda|$  is small enough
  and $  \lambda \neq 0, $  then $T-\lambda I$ is a Fredholm operator.
From \cite[Theorem 4.7]{GR} it follows that $T$ has a finite
essential ascent and finite essential descent.  Thus $T$ is a
B-Fredholm operator and $a_e(T)= d_e(T) \leq n.$ As $T$ is a
B-Fredholm operator, then from Theorem \ref{order}, we have  $ n
\leq a_e(T)= d_e(T)$  and so, $a_e(T)= d_e(T) = n.$

\ep \vspace{5mm}

 {\bf \underline{Application}}

 We give now an application of the previous results for the study of
the mean convergence in the Calkin algebra. For the uniform
ergodic theorem, we refer the reader to \cite[Theorem 1.5]{BUR}
and the references cited there. Here, using Theorem \ref{lifting},
we obtain easily  a general characterization of the convergence of
the sequence  $ (\Pi(M_n(T))_n$ in the Calkin algebra.

\bet \label{Pi-mean} Let $T \in L(X)$ and let $ M_n(T)= \frac{1+ T
+ T^2+...T^n}{n}, n \in \mathbb{N^*}.$ Then following conditions
are equivalent:

1- The sequence $ (\Pi(M_n(T))_n$ converges in the Calkin algebra.

2-  $\frac{ \mid\mid \Pi(T)^n\mid\mid}{n} \rightarrow 0$ as $
n\rightarrow \infty$ and there exists a power compact operator $K$
such that $ I-T+K$ is a B-Fredholm  operator, and the commutator
$[T,\,K]$ is compact.

3-  $\frac{ \mid\mid \Pi(T)^n\mid\mid}{n} \rightarrow 0$ as $
n\rightarrow \infty$ and there exists a power compact operator $K$
such that $ a_e(I-T+K)$ and $ d_e(I-T+K)$ are both finite  and the
commutator $[T,\,K]$ is compact.

4- $\frac{ \mid\mid \Pi(T)^n\mid\mid}{n} \rightarrow 0$ as $
n\rightarrow \infty$ and $ 1$ is a pole of the resolvent of
$\Pi(T).$

\eet

\bp $1)\Rightarrow 2)$  Assume that the sequence $ (\Pi(M_n(T))_n$
converges in the Calkin algebra. Then from \cite[Theorem
1.5]{BUR}, it follows that $\frac{ \mid\mid \Pi(T)^n\mid\mid}{n}
\rightarrow 0$ as $ n\rightarrow \infty$ and $ 0 $ is a pole of
the resolvent of $ \Pi(I-T).$ From Theorem \ref{lifting}, there
exists a power compact operator $K$ such that $ I-T+K$ is a
B-Fredholm operator and the commutator $[T,\,K]$ is compact.

$2) \Rightarrow 3)$  Assume $\frac{ \mid\mid \Pi(T)^n\mid\mid}{n}
\rightarrow 0$ as $ n\rightarrow \infty$ and there exists a power
compact operator $K$ such that $ I-T+K$ is a B-Fredholm operator.
Then $ a_e(I-T+K)$ and $ d_e(I-T+K)$ are both finite and the
commutator $[T,\,K]$ is compact.

$3) \Rightarrow 4)$ Assume  that $\frac{ \mid\mid
\Pi(T)^n\mid\mid}{n} \rightarrow 0$ as $ n\rightarrow \infty$ and
there exists a power compact operator $K$ such that $ a_e(I-T+K)$
and $ d_e(I-T+K)$ are both finite, and the commutator $[T,\,K]$ is
compact. Then  from  Theorem \ref{lifting},   $1$ is a pole of the
resolvent of $ \Pi(T).$

$4) \Rightarrow 1)$   Assume that  $\frac{ \mid\mid
\Pi(T)^n\mid\mid}{n} \rightarrow 0$ as $ n\rightarrow \infty$ and
 $1$ is a pole of the resolvent of $ \Pi(T).$ Using \cite[Theorem
1.5]{BUR}, it follows that  the sequence $ (\Pi(M_n(T))_n $
converges in the Calkin algebra.\\

\ep

Note that from \cite[Theorem 1.5]{BUR}, if $T$ satisfies one of
the conditions of Theorem \ref{Pi-mean}, then  $ 1$ is a pole of
the resolvent of $\Pi(T)$  of order less or equal  to $1.$

\bremark In the case of a Hilbert space, and from Corollary
\ref{Hilbert-case}, the operator $K$ of Theorem \ref{Pi-mean} can
be chosen to be compact. \eremark

 \baselineskip=12pt
\bigskip
\vspace{-5 mm }
 \baselineskip=12pt
\bigskip

{\small
\noindent Mohammed Berkani,\\
 \noindent Science faculty of Oujda,\\
\noindent University Mohammed I,\\
\noindent Laboratory LAGA, \\
\noindent 60000 Oujda, Morocco\\
\noindent berkanimo@aim.com,\\

\noindent Sne\v zana  \v{C}. \v{Z}ivkovi\'{c}-Zlatanovi\'{c}

\noindent University of Ni\v{s}, \\
\noindent Faculty of Sciences and Mathematics,\\
\noindent 18000 Ni\v{s}, Serbia\\
\noindent mladvlad@mts.rs

\end{document}